\documentclass[11pt]{article}
\usepackage{epsfig,amsmath,amssymb,stmaryrd,enumerate}
\usepackage{graphicx}

\usepackage{latexsym}
\usepackage{amsthm}
\usepackage{amsmath}
\usepackage{amsfonts, epsfig, amsmath, amssymb, color, amscd}
\usepackage[cmtip,arrow]{xy}
\usepackage{pb-diagram, pb-xy}
\usepackage{amssymb,epsfig}
\usepackage{color}
\usepackage[cmtip,arrow]{xy}
\usepackage{pb-diagram, pb-xy}
\usepackage{amssymb,epsfig,amsfonts}

\newfont{\sdbl}{msbm9}
\newfont{\dbl}{msbm10 at 12pt}
\theoremstyle{definition}

\newcommand{\ca}{{{\cal A}}}

\theoremstyle{plain}

\newtheorem{theo}{Theorem}[section]
\newtheorem{lemma}{Lemma}[section]

\newcommand{\veps}{\varepsilon}

\newcommand{\cs}{{{\cal S}}}

\newcommand{\cb}{{{\cal B}}}

\begin{document}

\title{Angular Billiard and Algebraic Birkhoff conjecture}

\author{{Misha Bialy\thanks{School of Mathematical Sciences, Raymond and Beverly Sackler Faculty of Exact Sciences, Tel Aviv University,
Israel},\quad Andrey E. Mironov\thanks{Sobolev Institute of Mathematics, Novosibirsk, Russia}
}}
\date{}
\maketitle

\begin{abstract}
In this paper we introduce a new dynamical system which we call
Angular billiard. It acts on the exterior points of a convex curve
in  Euclidean plane. In a neighborhood of the boundary curve this
system turns out to be dual to the Birkhoff billiard. Using this
system we get new results on algebraic Birkhoff conjecture on
integrable billiards.

\end{abstract}

\section{Introduction and Main Results}
\subsection{Motivation}
Birkhoof conjecture for convex billiards asks if there exist
integrable convex billiards other than ellipses. This paper is
motivated by the following two remarkable results on algebraic
version of this conjecture. The first one \cite{B}, due to S.V.
Bolotin (1990), deals with polynomial integrals of classical
Birkhoff billiards. The second one \cite{T2}, by S.L. Tabachnikov
(2008), is on polynomial integrals for the so called Outer
billiards. The methods of the proof of these results are different,
but both are algebro--geometric. We are particularly inspired by the
method of Tabachnikov and explain how to extend it to classical
Birkhoff billiards.

Our main tool is to introduce a new, the so called Angular billiard
which turns out to be dual to Birkhoff billiard in a neighborhood of
the boundary curve. Using Angular billiards we were able to extend
both of the results.

It is important to notice that our method is based on infinitesimal
expansion near the boundary. First, using Angular billiard we write
an equation depending on a small parameter $\varepsilon$ and then
extract from it the first non-trivial term. This gives us, see
Theorem \ref{remarkable}, a remarkable identity which is further
studied in algebro--geometric terms. It is plausible,  that the
analysis of the higher order terms in the expansion (they become
much more complicated) may lead to the complete proof of the
algebraic Birkhoff conjecture. It is worth mentioning recent results
\cite{Ka1}, \cite{Ka2} around the Birkhoff conjecture where certain
expansions near the boundary were used also.

Let us remark also that since our method is infinitesimal it is
related to the question when the family of caustics near the
boundary of integrable Birkhoff billiard is algebraic.

Let us fix the notations. Let $\gamma$ be a simple, closed,
piecewise smooth  curve in the Euclidean plane bounding the domain
$\Omega.$ Given a tangent vector $v\in T_q\Omega$, Birkhoff billiard
flow $g^t$ transports $v$ parallel along straight line in the
direction of $v$ till the collision with the boundary, where the
vector gets reflected according to the low of geometric optics:
$$v_+=v_--2n <v_- ,n>,$$
where $n$ is an outward normal vector to $\gamma$ and $<\cdot\
,\cdot>$ is the scalar product. Obviously, the energy, that is
$|v|^2,$ is conserved by dynamics and it is a very classical
question by Birkhoff when there exists a non constant additional
function $ \Phi=\Phi(q,v)$ defined on the unite energy level
$\{|v|=1\}$ which is \textit{polynomial} in the components of vector
$v=(v_x, v_y)$ such that it is a first integral of the billiard flow
$g^t,$ i.e. $\Phi\circ g^t=\Phi.$ We refer to \cite{KT}, \cite{T2},
and also to Section 4 below for more details.

\vskip6mm \noindent $\mathbf{Theorem}$ (Bolotin, 1990) {\it Assume
that Birkhoff billiard inside $\gamma$ admits a non-constant
polynomial integral $\Phi$ on the energy level $\{|v|=1\}.$ It then
follows that all smooth arcs $\gamma_i$ of $\gamma$ are necessarily
algebraic curves. Moreover, let $\tilde{\gamma}_i$
 be the irreducible component of $\gamma_i$ in $\mathbb{C}P^2.$ Then, the
following alternative holds: either every $\tilde{\gamma}_i$ is of
degree 1 or 2, or every $\tilde{\gamma}_i$ necessarily contains
singular points. In particular, if $\gamma$ contains only one smooth
closed arc then it is either an ellipse, or the irreducible
component, $\tilde{\gamma},$ of $\gamma$ in $\mathbb{C}P^2$ is a
singular curve.} \vskip6mm

It turns out that the Angular billiards introduced below allow us to
extend this result in terms of the dual curve $\Gamma.$ Throughout
the paper we fix Euclidean coordinates $(x,y)$ on the Euclidean
plane centered at a fixed point $O$ inside $\gamma,$ and extend them
to homogeneous coordinates $(x:y:z)$ in $\mathbb{C}P^2.$ We shall
also identify $\mathbb{C}P^2$ and the dual projective space with the
help of the quadric $\{x^2+y^2=z^2\}$ corresponding to Euclidean structure.
Then the polar duality in the
Euclidean plane coincides with the duality in $\mathbb{C}P^2$ (see
\cite{F} and also below). With these remarks our main result reads:
\vskip6mm

\noindent{\bf Theorem 1} {\it Suppose that Birkhoff billiard inside
$\gamma$ admits a non-constant polynomial integral $\Phi$ on the
energy level $\{|v|=1\}.$ For every smooth non-linear piece $\gamma_i$
 of $\gamma,$ let $\Gamma_i$ be the dual curve.
Denote $\tilde{\Gamma}_i $ the irreducible
component in $\mathbb{C}P^2$ of $\Gamma_i.$ Then,
either every $\tilde{\Gamma}_i$ has degree 2, or every
$\tilde{\Gamma}_i$ necessarily contains singular points. Moreover,
all singular and inflection points of $\tilde{\Gamma}_i$ in
$\mathbb{C}P^2$ belong to the union of two lines defined by the equations
$$
L_1=\{ x+iy=0\}, \qquad L_2=\{x-iy=0\}.
$$}
%\vskip6mm \noindent{\bf Remark} We explain

\vskip6mm \noindent{\bf Corollary} $\textbf{1.}$ {\it If the
Birkhoff billiard inside $\gamma$ is integrable  with an integral which is polynomial in
$v$, then $\tilde{\gamma}$ does not have two real algebraic
ovals having a common tangent line.}

\vskip6mm Indeed, if there were
two such ovals $\gamma$ and $\gamma_1$ (see Fig. 1), then the point
$(x,y)$ on $\Gamma$ dual to the common tangent line $\tau$
would be a real singular point of $\Gamma$ different from $O,$ and hence
$x^2+y^2\ne 0.$  It then follows from Theorem 1 that the Birkhoff billiard
inside $\gamma$ does not admit polynomial integral on the energy
level $|v|=1.$

\vskip20mm

\begin{picture}(170,100)(-100,-50)

\qbezier(40,50)(41,86)(95,87) \qbezier(95,87)(149,86)(150,50)
\qbezier(150,50)(149,15)(95,14) \qbezier(95,14)(41,15)(40,50)
\put(165,32){\shortstack{$\gamma_1$}}
\put(33,30){\shortstack{$\gamma$}}

\put(170,49){\circle{20}} \put(110,91){\line(2,-1){83}}
\put(110,92){\shortstack{$\tau$}}
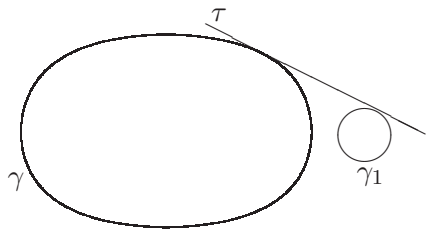
\put(10,-5){\shortstack{Fig. 1 {\footnotesize Non-integrable
Birkhoff billiard inside $\gamma.$}}}
\end{picture}

Let us consider the simplest example to Corollary 1. We take real algebraic curve
$$
 y^2=F(x)=(x-x_1)(x-x_2)(x-x_3)(x-x_4)f(x),\quad x_j\in{\mathbb R},
 \quad x_1<x_2<x_3<x_4,
$$
where $f(x)$ is a real polynomial such that $F(x)> 0$ for $x\in
(x_1,x_2)$ and $x\in (x_3,x_4).$
%($f(x)$ can have multiple roots).
Then the Birkhoff billiard inside the real analytic oval
$$
 \gamma=\{(x,\pm\sqrt{F(x)}),\ x\in[x_1,x_2]\}
$$
does not admit polynomial integral, since the algebraic curve has
another analytic oval $\gamma_1$ as in Fig. 1
$$
 \gamma_1=\{(x,\pm\sqrt{F(x)}),\ x\in[x_3,x_4]\}.
$$
Moreover, since the algebraic curve $\tilde{\gamma}\subset{\mathbb
C}P^2$ is always singular, then Theorem by Bolotin does not apply.
Notice, that if the polynomial $f$ has multiple roots then
singularities of $\tilde{\gamma}$ are complicated.

\vskip6mm \noindent{\bf Corollary}$\textbf{2.}$ {\it  Assume that
$\tilde{\Gamma}$ is a non-singular curve (of degree $>2$) in
$\mathbb{C}P^2$ and has a  smooth real oval $\Gamma$ (for example,
$\tilde{\Gamma}$ is a nonsingular cubic). Then the dual curve
$\gamma$ is also an oval and Birkhoff billiard inside $\gamma$ is
not integrable by Theorem 1. Notice, that in this case Bolotin's
theorem does not apply, since $\tilde{\gamma}$ is necessarily
singular in this case. The inflection points of $\tilde{\Gamma}$
correspond to singular points of $\tilde{\gamma}.$ }

\vskip6mm \noindent{\bf Corollary}$\textbf{3.}$ {\it Let  $\gamma$
be the dual of Fermat oval $\Gamma=\{x^{2n}+y^{2n}=1,\ n>1 \}.$
Notice that $\tilde{\Gamma}$ is irreducible, non-singular curve and so by
Theorem 1 the Birkhoff billiard inside $\gamma$ is not integrable.
One can easily compute that in this case the oval $\gamma$
 can be written as follows:
 $$
 \gamma=\{x^{\frac{2n}{2n-1}}+ y^{\frac{2n}{2n-1}}=1\}.
 $$
Therefore (the algebraic curve) $\gamma$ is a strictly convex $C^1$
curve in the plane which has 4 singular points $(\pm 1,0),(0,\pm 1)$
corresponding to 4 inflection points $(\pm 1,0),(0,\pm 1)$ of
$\Gamma.$ So Bolotin's theorem does not apply in this case also. }

\vskip6mm

 In the sequel we shall assume, for simplicity of the exposition,
 that $\gamma$ has only one smooth piece.
 Proof of the general case is completely analogous. Proof of Theorem 1
 is contained in Section 6.

%%%%%%%%%%%%%%%%%%%%%%%%%%%%%%%%%%%%%%%%%%%%%%%%%%%
\subsection{Introducing Angular billiards} Let $D\subset{\mathbb R}^2$ be a convex
domain bounded by a smooth curve $\Gamma,$  $O$ be a fixed point
inside $\Gamma.$ We shall always choose a counter clockwise
orientation on $\Gamma.$ Let us take an arbitrary point $A$ outside
$\Gamma,$ $A\in U={\mathbb R}^2\backslash D.$  There are two tangent
lines from $A$ to $\Gamma.$ We choose the right one (if one looks at
$\Gamma$ from $A$), which agrees with orientation of $\Gamma$. Let
$T$ be the tangency point.

\vskip23mm

\begin{picture}(170,100)(-100,-50)

\qbezier(40,50)(41,86)(95,87)
\qbezier(95,87)(149,86)(150,50)
\qbezier(150,50)(149,15)(95,14)
\qbezier(95,14)(41,15)(40,50)

\put(40,50){\vector(0,-2){2}}

\put(147,25){\shortstack{$\Gamma$}}

\put(76,62){\shortstack{$l_A$}}

\qbezier(128,65)(121,65)(118,59)

\put(115,69){\shortstack{$\alpha$}}

\put(106,64){\shortstack{$\alpha$}}

\put(130,55){\circle*{2}}
\put(125,43){\shortstack{$O$}}

\put(124,87){\line(-4,0){89}}

\put(124,87){\circle*{2}}

\put(119,90){\shortstack{$A$}}

\put(130,55){\line(-1,5){6.5}}

\put(130,55){\line(-1,1){32}}

\put(98,87){\circle*{2}}
\put(94,90){\shortstack{$T$}}

\put(130,55){\line(-3,1){94}}

\put(35,87){\circle*{2}}
\put(30,90){\shortstack{$B$}}

\put(15,-5){\shortstack{Fig. 2 {\footnotesize $\angle AOT=\angle
TOB,\quad \ca (A)=B.$}}}

\end{picture}

If $OA$ is not orthogonal to $OT,$ then there is unique line $l_A$
(different from $OA$) such that the angle between $OA$ and $OT$ is
equal to the angle between lines $OT$ and $l_A.$ If $l_A$ is not
parallel to $AT,$ then there is the intersection point $B$ of the
tangent line and $l_A$ (see Fig. 2). Let ${\cal S}\subset U$ be a
curve defined by the condition
$$
 {\cal S}=\{A: l_A\| AT\}.
$$
We introduce the mapping
$$
 \ca:U\backslash {\cal S}\rightarrow U,\qquad \ca(A)=B.
$$
Here we assume that if $OA$ is orthogonal to $OT$ (i.e.
$OA=l_A$), then $\ca(A)=A.$ On the curve ${\cal S}$ the mapping
$\ca$ can not be defined. We call the mapping $\ca$ the {\it Angular
Billiard map of $\Gamma.$}

We shall consider the mapping $\ca$ in details in the next section
2.

We shall explain now, that on the Euclidean plane, in a neighborhood
of the boundary curve, Birkhoff and Angular billiards are dual to
each other. In particular, if $\Gamma$ is sufficiently smooth
strictly convex curve, then KAM theory implies that in a
neighborhood of $\Gamma$ there are infinitely many invariant curves.
Moreover, we shall prove in Section 3 that near the boundary
$\Gamma$ Angular billiard map is a twist symplectic map of the
cylinder.

\vskip6mm \noindent $\mathbf{Remark.}$ {\it The reader shouldn't
confuse with Outer billiards studied by S. Tabachnikov in
(\cite{T1}) where they were called "Dual" because they indeed become
dual on the sphere. However, on the Euclidean plane, in a
neighborhood of the boundary, Dual to Birkhoff billiard is in fact
the Angular billiard. } \vskip6mm

We shall discuss KAM result and establish twist property for Angular
billiards in details in section 3.

\subsection{Duality in a neighborhood of the boundary}
Birkhoff billiard map $\cb$ acts (by the reflection in $\gamma$) on
the space of \emph{oriented} straight lines intersecting $\gamma.$
Recall, we have fixed the point $O$ inside $\gamma$ and Euclidean
coordinates $(x,y)$ centered at $O.$ Let us call an oriented line
not passing through $O$ positive (or respectively negative) if the
momentum $\sigma$ with respect to $O$ of any positively oriented
vector $v$ tangent to the line is positive (or negative
respectively):
$$\sigma(v)=xv_y-yv_x>0\ (<0).$$
Given an oriented line which passes sufficiently close to the
boundary $\gamma$ it's sign remains unchanged under the iterations
of the Birkhoff billiard. (For other lines not close to the boundary
the sign of course may change as we discuss below.)

Let us remind the correspondence of polar duality with respect to
the unite circle. It acts on \emph{non-oriented} lines. Given a line
not passing through $O$ we write $l$ in the form \newline $<n,x>=p,\
p>0,$ where $n$ is an outward unite normal to $l$ (that is $(n,v)$
is a positive basis, $v$ is the unite tangent vector to $l$ with
positive momentum) then the dual point corresponding to $l$ is by
definition $L=n/p.$ In other words, the corresponding point $L$ lies
on the normal radius to the line at the distance equal $1/p$:
\begin{equation}
\label{duality}l= \{<n,x>=p,\ p>0\}\leftrightarrow L=n/p.
\end{equation}
 As usual, we shall denote everywhere by
small letters the lines and by capital letters the corresponding
dual points.

We denote by $\Gamma$ the dual curve to $\gamma$ consisting of
points which are dual to the tangent lines of $\gamma.$ The origin
$O$ remains inside $\Gamma.$ Remarkably, the duality preserves the
incidence relation and dual to $\Gamma$ is $\gamma$ again. More
precisely, if $l$ is tangent to $\gamma$ at $Q$ then the dual line
$q$ is tangent to $\Gamma$ at $L$ (see Fig. 3).

Furthermore, if $a, b$ are two oriented positive lines in the plane
so that the Birkhoff billiard map $\cb$ transforms $a$ to $b.$ Let
$Q\in\gamma$ be the point of reflection and let $l$ be the tangent
line to $\gamma$ at the reflection point $Q.$ Then the dual points
$A,B$ lie on the line $q$ which is tangent to $\Gamma$ at $L.$
Moreover, we claim that the angles $AOL$ and $BOL$ are equal
implying the Angular billiard rule:
$$\ca(A)=B.$$

Indeed, the angle $\alpha$ between the lines $a$ and $l$ equals the
angle $\angle AOL$ because $OA$ is normal to $a$ and $OL$ is normal
to $l,$ by definition of polar duality. Analogously, the angle
$\beta$ between the lines $l$ and $b$ equals the angle $\angle LOB$
(see Fig. 3). Moreover, $\alpha$ and $\beta$ are equal by the
reflection law of Birkhoff billiard. This proves the claim.

\vskip35mm

\begin{picture}(170,100)(-100,-50)

\qbezier(-50,50)(-49,86)(5,87)
\qbezier(5,87)(59,86)(60,50)

\qbezier(-3,95)(-8,87)(-3,79)
\put(-16,90){\shortstack{$\alpha$}}
\put(-16,74){\shortstack{$\beta$}}

\put(5,87){\circle*{2}}
\put(0,90){\shortstack{$Q$}}
\put(42,50){\vector(-1,1){37}}
\put(5,87){\vector(-1,1){37}}
\put(5,87){\vector(-1,-1){37}}
\put(25,67){\shortstack{$a$}}
\put(-14,107){\shortstack{$a$}}

\put(-40,50){\shortstack{$b$}}

\put(60,87){\line(-4,0){110}}
\put(55,90){\shortstack{$l$}}
\put(62,55){\shortstack{$\gamma$}}

\put(5,50){\circle*{2}}
\put(-5,40){\shortstack{$O$}}

\qbezier(177,58)(185,63)(193,58)
\put(187,63){\shortstack{$\alpha$}}
\put(175,63){\shortstack{$\beta$}}

\qbezier(130,50)(131,86)(185,87)
\qbezier(185,87)(239,86)(240,50)
\put(185,87){\circle*{2}}
\put(180,90){\shortstack{$L$}}
\put(240,87){\line(-4,0){110}}
\put(232,92){\shortstack{$q$}}
\put(222,87){\circle*{2}}
\put(148,87){\circle*{2}}
\put(215,90){\shortstack{$A$}}
\put(145,90){\shortstack{$B$}}
\put(242,55){\shortstack{$\Gamma$}}

\put(185,50){\circle*{2}}
\put(175,40){\shortstack{$O$}}

\put(185,50){\line(0,3){38}}
\put(185,50){\line(-1,1){37}}
\put(185,50){\line(1,1){37}}

\put(55,-5){\shortstack{Fig. 3 {\footnotesize Polar duality; $\beta=\alpha.$}}}

\end{picture}

%%%%%%%%%%%%%%%%%%%%%%%%%%%%%%%%%%%%%%%%%%%%%%%%%%%%%%%%%%%%%%%
In the following section we consider dynamics of $\ca$ in more
details. We can summarize the relation of $\ca$ and $\cb$ as
follows:

\vskip8mm

\noindent {\bf Theorem 2}  {\it Suppose $\cb$ acts on an oriented
line $l_0:$ $\cb(l_0)=l_1,\ \cb^{-1}(l_0)=l_{-1}.$   Then the
Angular billiard acts on the dual points as follows: $$
\ca(L_0)=L_1,\quad if\ l_0\ is\ positive ,\quad\ca(L_0)=L_{-1},\quad
if\  l_0\  is\ negative.$$ In terms of the Angular billiard this
reads as follows. Let $A$ be a point outside $\Gamma,$ $B=\ca(A).$
Let $a,b$ are positively oriented dual lines and $l$ is the positive
tangent line to $\Gamma$ passing through $A$ with the tangency point
$T.$ Then
$$
\cb (a)=b,\ if\ T\in[A;B], \quad and\quad \cb(a)=-b,\ if\ T\notin
[A;B].$$}

The illustration of Theorem 2 is contained in Section 2.

\vskip6mm
\subsection{Integrability of Angular billiard}

We shall call the Angular billiard {\it integrable} if there is a
function
$$
 G:U\backslash {\cal S} \rightarrow {\mathbb R}
$$
such that
$$
G(A)=G(\ca (A)),\qquad\forall A\in U\backslash {\cal S}.
$$

\vskip8mm

\noindent {\bf Example 1} {\it Let $\Gamma$ be an ellipse defined
by the equation
$$
 \frac{x^2}{a^2}+\frac{y^2}{b^2}=1,
$$
and let $O(x_0,y_0)$ be an arbitrary point inside $\Gamma.$ Then the Angular billiard is integrable and the integral has the form
$$
 G(x,y)=\frac{F(x,y)}{(x-x_0)^2+(y-y_0)^2},
 \qquad F(x,y)=\frac{x^2}{a^2}+\frac{y^2}{b^2}-1.
$$}

\vskip6mm

It is a remarkable fact that the integrability of Angular billiard for the ellipse
does not depend on the choice of the point $O.$

Let us suppose now that the Birkhoff billiard flow admits a
polynomial integral $\Phi$ on the energy level $\{|v|=1\}.$ Then it
is shown in \cite{B}, \cite{KT}, that one can assume that
$\Phi(q,v)$ is a \emph{homogeneous} polynomial of a certain
\emph{even} degree $n$ in $\sigma(v)=xv_y-yv_x, \ v_x,v_y:$
$$\Phi=\Phi(\sigma,v_x,v_y).
$$  In order not to overload the notations we denote them by the
same letter $\Phi.$ Moreover, one can modify $\Phi$ in such a way
that $\Phi$ vanishes on any tangent vector to the boundary $\gamma.$
See section 4 for details.

\vskip8mm

\noindent{\bf Theorem 3} {\it Let $\gamma$ be a convex closed curve
and $\Phi(\sigma,v_x,v_y)$ be a homogeneous polynomial integral of
even degree $n,$ vanishing on the tangent vectors to the boundary
$\gamma.$ Then the Angular billiard corresponding to $\Gamma$ is
also integrable with the integral of the form
$$
 G(x,y)=\frac{F(x,y)}{(\sqrt{x^2+y^2})^n},\qquad F(x,y)=\Phi(1,-y,x),
$$
where $F$ is a (non-homogeneous) polynomial of degree $n.$ Moreover
$F$ vanishes on $\Gamma.$}

\vskip8mm

\noindent {\bf Example 2} {\it Let $\gamma$ be an ellipse defined by the equation
$$
 a^2x^2+b^2y^2-1=0,\quad a\geq b.
$$
The Birkhoff billiard is integrable inside $\gamma,$ and the first
integral satisfying the condition that it vanishes on tangent
vectors to $\gamma,$ has the form:
$$
 \Phi(\sigma,v)=\frac{v_x^2}{b^2}+\frac{v_y^2}{a^2}-\sigma^2=
 \frac{v_x^2}{b^2}+\frac{v_y^2}{a^2}-(xv_y-yv_x)^2.
$$
The dual curve $\Gamma$ to the ellipse $\gamma$ is given by the
equation
$$
 \frac{x^2}{a^2}+\frac{y^2}{b^2}-1=0.
$$
According to Theorem 3 the Angular billiard for $\Gamma$ has the first integral
$$
 G(x,y)=\frac{\frac{x^2}{a^2}+\frac{y^2}{b^2}-1}{x^2+y^2}.
$$
This is exactly the integral which is given in Example 1.}

\vskip8mm

Let us mention that Theorem 1 can be  effectively improved when the
integral has small degree. More precisely, we have proved the
following

\vskip6mm

\textbf{Theorem}
 {\it For a smooth closed curve $\gamma$ different from ellipse
Birkhoff billiard inside $\gamma$ does not admit polynomial integral
of degree
 4.}
 \vskip6mm
It is worth mentioning  the piecewise smooth examples of
\cite{greki}, where the curve $\gamma$ is composed of confocal arcs
and straight segments, where there is an integral of degree four.
Much more general questions and phenomena were later on studied, see
\cite{dragovich}.

We shall publish the proof of this theorem later on.

 \vskip6mm

The plan of the paper is as follows: In Section 2 we study the
Angular billiard mapping $\ca$ and we explain the connection of
Theorem 2 between Angular billiards and Birkhoff billiards. In
section 3 we prove KAM and twist property for angular billiards. In
section 4 we remind facts on Polynomial integrals. In section 5
Theorem 3 is proved. In Section 6 we deduce a remarkable equation
(Theorem \ref{remarkable}) and prove Theorem 1.

\section*{Acknowledgements}
It is a pleasure to thank Eugene Shustin for useful consultations.
M. Bialy was supported in part by ISF grant 162/15, A.E. Mironov was
supported by RSF (grant 14-11-00441).

Essential part of the paper was obtained during A.E. Mironov visit to
Tel Aviv University, it is a pleasure to thank TAU for excellent
working conditions.

\section{Angular Billiard Map}

Theorem 2 follows from the following geometric analysis of the
action of the Angular billiard map $\ca.$ First of all let us study
the action  of $\ca$ on points which lie on a tangent line to
$\Gamma.$ Choose a point $A\in U$ and the positive tangent line to
$\Gamma.$ There are three possibilities (see Fig. 2): \vskip5mm 1.
In the first case, $\angle ATO=\frac{\pi}{2}$ the mapping $\ca$ is a
bijection between half lines $(-\infty;T)$ and $(T;+\infty).$
\vskip8mm
 2. In the second case $\angle ATO>\frac{\pi}{2}$. Then
there is a point $M$ on the tangent line such that $\angle
M_1OT=\angle MOT,$ where $OM_1\| TA$ (see Fig. 4). The mapping $\ca$
maps the half  line $(-\infty;T)$ on the interval $(T;M).$
\newpage
\vskip35mm

\begin{picture}(170,100)(-100,-50)
\qbezier(40,50)(41,86)(95,87)
\qbezier(95,87)(149,86)(150,50)
\qbezier(150,50)(149,15)(95,14)
\qbezier(95,14)(41,15)(40,50)

\put(40,50){\vector(0,-2){2}}

\put(147,65){\shortstack{$\Gamma$}}

\put(128,55){\circle*{2}}
\put(128,57){\shortstack{$O$}}

\put(128,55){\line(-4,0){63}}
\put(78,55){\circle*{2}}
\put(71,60){\shortstack{$M_1$}}

\put(95,14){\circle*{2}}
\put(90,1){\shortstack{$T$}}

\put(128,55){\line(-4,-5){33}}

\put(128,55){\line(-5,-3){69}}

\put(60,14){\circle*{2}}
\put(53,1){\shortstack{$A$}}

\put(21,14){\line(4,0){135}}
\put(20,3){\shortstack{$-\infty$}}

\put(128,55){\line(-1,-3){14}}
\put(114,14){\circle*{2}}
\put(108,1){\shortstack{$B$}}

\put(128,55){\line(1,-3){14}}
\put(141.5,14){\circle*{2}}
\put(133,1){\shortstack{$M$}}

%\put(68,32){\shortstack{$\alpha$}}
%\put(76,36){\shortstack{$\alpha$}}

\qbezier(113,46)(117,40)(124,42)
\put(111,32){\shortstack{$\alpha$}}
\put(103,36){\shortstack{$\alpha$}}

\put(5,-20){\shortstack{Fig. 4 {\footnotesize $\angle
ATO>\frac{\pi}{2},\quad \ca (A)=B,\quad \ca (-\infty; T)=(T;M).$}}}

\end{picture}

3. In the third case $\angle ATO<\frac{\pi}{2}.$ This case is more interesting. On the tangent line there is a
unique point $S$ such that $\angle SOT=\angle S_1OT,$ where $OS_1\|
ST.$ Then for arbitrary point $A$ from the interval $(S;T)$ there is
a unique point $B$ from the another half of the tangent line such
that
$$\angle AOT=\angle TOB.$$
The mapping $\ca$ maps the interval $ST$ on the half tangent line
$(T;+\infty)$ containing $B$ (see Fig.~5).

\vskip20mm

\begin{picture}(170,100)(-100,-50)

\qbezier(40,50)(41,86)(95,87)
\qbezier(95,87)(149,86)(150,50)
\qbezier(150,50)(149,15)(95,14)
\qbezier(95,14)(41,15)(40,50)

\put(40,50){\vector(0,-2){2}}

\put(147,25){\shortstack{$\Gamma$}}

\qbezier(128,65)(121,65)(118,59)

\put(116,69){\shortstack{$\alpha$}}
\put(106,64){\shortstack{$\alpha$}}

\put(130,55){\circle*{2}}
\put(125,43){\shortstack{$O$}}
\put(150,87){\line(-4,0){155}}

\put(130,55){\line(-4,0){65}}
\put(85,55){\circle*{2}}
\put(79,43){\shortstack{$S_1$}}

\put(-5,90){\shortstack{$+\infty$}}

\put(124,87){\circle*{2}}
\put(118,90){\shortstack{$A$}}
\put(130,55){\line(-1,5){6.5}}

\put(130,55){\line(1,5){6.5}}
\put(137,87){\circle*{2}}
\put(132,90){\shortstack{$S$}}

\put(130,55){\line(-1,1){32}}
\put(98,87){\circle*{2}}
\put(93,90){\shortstack{$T$}}
\put(130,55){\line(-3,1){94}}
\put(35,87){\circle*{2}}
\put(30,90){\shortstack{$B$}}

\put(15,-5){\shortstack{Fig. 5 {\footnotesize $\angle
ATO<\frac{\pi}{2},\quad \ca (A)=B,\quad \ca(S;T)=(\infty; T).$}}}

\end{picture}

Furthermore, in the third case there exists a point $P$ on the half
tangent line such that $OP\perp OT$ (see Fig. 6).

For arbitrary point $A$ from the interval $(-\infty ;P)$ there is a
unique point $B$ from the interval $(P;S)$ such that the angle
between lines $AO$ and $OT$ is equal to the angle between lines $OT$
and $OB$ or since $OP\perp OT$ it is equivalent to $\angle
AOP=\angle BOP.$

Hence,
$$
 \ca (A)=B,\quad \ca (B)=A,\quad \ca(P)=P,
$$
and $\ca$ maps the half line $(-\infty; P)$ on the interval $(S;P).$
When the tangent line varies, the points $P$ form a curve in $U$
$$
 {\cal P}=\{P\in U| OP\perp OT\}\subset U.
$$

\vskip25mm
\begin{picture}(170,100)(-100,-50)

\qbezier(40,50)(41,86)(95,87)
\qbezier(95,87)(149,86)(150,50)
\qbezier(150,50)(149,15)(95,14)
\qbezier(95,14)(41,15)(40,50)

\put(40,50){\vector(0,-2){2}}

\put(147,25){\shortstack{$\Gamma$}}

\put(130,55){\circle*{2}}
\put(125,43){\shortstack{$O$}}
\put(195,87){\line(-4,0){135}}

\put(130,55){\line(-4,0){65}}
\put(85,55){\circle*{2}}
\put(79,43){\shortstack{$S_1$}}

\put(130,55){\line(1,6){5.5}}
\put(135,87){\circle*{2}}
\put(131,90){\shortstack{$S$}}

\put(130,55){\line(4,5){25}}
\put(155,87){\circle*{2}}
\put(149,90){\shortstack{$P$}}

\put(130,55){\line(1,2){16}}
\put(146,87){\circle*{2}}
\put(140,90){\shortstack{$B$}}

\put(130,55){\line(6,5){38}}
\put(168,87){\circle*{2}}
\put(163,90){\shortstack{$A$}}

\put(185,90){\shortstack{$-\infty$}}

\put(130,55){\line(-1,1){32}}
\put(98,87){\circle*{2}}
\put(93,90){\shortstack{$T$}}

\put(-10,-5){\shortstack{Fig. 6 {\footnotesize $\angle
ATO<\frac{\pi}{2},\ \ca (A)=B,\ \ca (B)=A, \ \ca(-\infty;
P)=(P;S).$}}}

\end{picture}

Thus, generally there is a domain $V\subset U$ such that the
restriction $\ca$ on $V$ is an involution (points $A,P$ and $B$
belong to $V$), $\ca^2=id$ and the curve ${\cal P}\subset V$ is
fixed under this involution.

Let us consider the simplest example.

\vskip10mm

\noindent {\bf Example 3} {\it Let us consider a circle of radius $r$
given by the equation $x^2+y^2=r^2$ and let $O$ be a point with
coordinates $(x_0,0),\ 0\leq x_0\leq r$ (see Fig. 7). Then the
curve ${\cal S}$ is a vertical line $x=\frac{r^2+x_0^2}{2x_0},$ and
the curve ${\cal P}$ is defined by
$$
 y(x)=\frac{(x-\frac{r^2}{x_0})(x-x_0)}{\sqrt{(x-\frac{r^2-rx_0+x_0^2}{x_0})(\frac{r^2+rx_0+x_0^2}{x_0}-x)}}.
$$
The curve ${\cal P}$ has two vertical asymptotic lines. Here
$$
 V=\{(x,y)\in{\mathbb R}^2|\ x>\frac{r^2+x_0^2}{2x_0}\}.
$$
}

\vskip25mm
\begin{picture}(170,100)(-100,-10)

\put(45,87){\circle{80}}
\put(15,70){\shortstack{$\Gamma$}}

\put(55,87){\circle*{2}}
\put(51,77){\shortstack{$O$}}

\put(90,45){\line(0,4){90}}
\put(92,50){\shortstack{${\cal S}$}}

\qbezier(160,135)(161,92)(140,87)
\qbezier(140,87)(120,82)(120,44)
\put(123,50){\shortstack{${\cal P}$}}

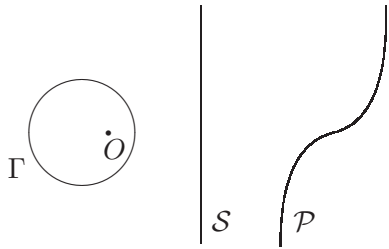
\put(20,15){\shortstack{Fig. 7 {\footnotesize Angular billiards on a circle}}}

\end{picture}

Let us explain the connection between Angular billiards and the
Birkhoff billiards. For this we construct a mapping $\tilde{\ca}$ on
the set of non-oriented lines with non-empty intersections with the domain
$\Omega$ as follows:

 Take an
arbitrary line $l_1$ not passing through $O,$ and a positive vector
$v_1(a_1,b_1)$ tangent to $l_1$ such that the angular momentum of
$v_1$ with respect to $O$ is positive
$$
 \sigma(v_1)=xb_1-ya_1> 0.
$$

\vskip20mm

\begin{picture}(170,100)(-100,-50)

\qbezier(40,50)(41,86)(95,87)
\qbezier(95,87)(149,86)(150,50)
\qbezier(150,50)(149,15)(95,14)
\qbezier(95,14)(41,15)(40,50)

\put(40,50){\vector(0,-2){2}}

\put(147,25){\shortstack{$\Gamma$}}

\put(60,55){\circle*{2}} \put(55,43){\shortstack{$O$}}

\put(120,84){\circle*{2}}
\put(117,90){\shortstack{$T$}}

\put(136,35){\vector(-1,3){16}} \put(136,36){\shortstack{$l_1$}}
\put(128,62){\shortstack{$v_1$}}

\put(91,41){\vector(2,3){29}}
\put(80,42){\shortstack{$l_2$}}
\put(90,62){\shortstack{$v_2$}}

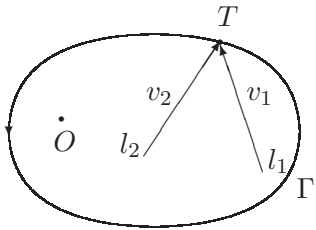
\put(5,-5){\shortstack{Fig. 8 {\footnotesize $\sigma(v_1)>0,\quad
\sigma(v_2)>0,\quad \tilde{\ca}(l_1)=l_2,\quad
\tilde{\ca}(l_2)=l_1.$}}}

\end{picture}

After the reflection we get another line $l_2$ and we choose
positive orientation on $l_2$ also. We set
$$\tilde{\ca}(l_1):=l_2.$$

 If the angle between $l_1$ and the tangent
line at the reflection point is rather small, we see that
$\tilde{\ca}=\cb$, so we obtain lines $l_n=\tilde{\ca}(l_{n-1})$ of
the trajectories of the usual Birkhoff billiard. If $l_1$ is
orthogonal to $\gamma$ at the reflection point, we have $
\tilde{\ca}(l_1)=l_1. $ If the reflected line $\cb(l_1)$ gets
negative orientation after the reflection (see Fig. 8), then we have
$$
 \tilde{\ca}(l_1)=l_2,\qquad \tilde{\ca}(l_2)=l_1.
$$

It then follows that  $\tilde{\ca}$ is transformed into $\ca$ with
the help of polar duality. In particular the curve ${\cal S}$ for
Angular billiards corresponds to the lines passing through $O,$ the curve
${\cal P}$ corresponds to the set of lines orthogonal to the
boundary $\gamma.$

\section{ Twist property and KAM for Angular billiard}

In this section we shall assume that $\gamma$ and therefore also
$\Gamma$ are smooth convex curves of everywhere positive curvatures.

It is known that for Birkhoff billiards there are two natural
choices of symplectic coordinates on the space of oriented lines in
the plane (see \cite{T2} for details). The first is $(\varphi,p)$
where $\varphi$ is the angle between the vector $(1,0)$  and the
positive normal to the line and $p=\sigma(v),$ where $v$ is the
unite vector in the direction of the line. In other words, $p$ is a
signed distance to the line from the origin.

The second choice of coordinates depends on a curve, in our case
$\gamma.$ The coordinates $(s,\cos\theta)$ where $s$ is a distance
traveled along the curve to the intersection point with the line and
$\theta$ is the angle between tangent vector to the curve with the
line. It is very useful fact that the billiard ball map is a twist
map with respect to the second choice of coordinates, the generating
function then is just the Euclidean Length (see \cite{T2}, \cite{katok} for twist maps
and generating functions). This fact is used in many results on
billiards. Theorem 5 below tells that near the boundary the billiard
map is a twist map also for the first choice of coordinates and
gives the generating function. It would be interesting to find new
variational applications of this new twist property.

Let us use the following notations. Let $r=r(\varphi)$ be
the equation of $\Gamma$ in the polar coordinates on the plane $\mathbb{R}^2(x,y).$ Here $r(\varphi)$ is
a position function of the curve $\Gamma,$ while $p(\varphi)=1/
r(\varphi)$ is the supporting function of $\gamma.$

\subsection{Twist condition and generating function}
We shall use complex notations

$$\Gamma(\varphi)=r(\varphi)e^{i\varphi},\quad
\dot{\Gamma}(\varphi)=(\dot{r}(\varphi)+ir(\varphi))e^{i\varphi}.$$

Let $M_1=r_1e^{i\varphi_1}, M_2=r_2e^{i\varphi_2}$ be two points so
that $\ca(M_1)=M_2.$ Let
$\bar{M}=r(\bar{\varphi})e^{i\bar{\varphi}},
\bar{\varphi}=\frac{\varphi_1+\varphi_2}{2}$ be the tangency point
of the segment $[M_1,M_2]$ with the curve $\Gamma.$ Recall $\ca$
denotes Angular billiard map.

We write for the points $M_1, M_2$
$$M_1=r_1e^{i\varphi_1}=\Gamma(\bar{\varphi})+
a\dot{\Gamma}(\bar{\varphi})=(r(\bar{\varphi})+a(\dot{r}(\bar{\varphi})
+ir(\bar{\varphi})))e^{i\bar{\varphi}},
$$
$$
M_2=r_2e^{i\varphi_2}=\Gamma(\bar{\varphi})+
b\dot{\Gamma}(\bar{\varphi})=
(r(\bar{\varphi})+b(\dot{r}(\bar{\varphi})
+ir(\bar{\varphi})))e^{i\bar{\varphi}}.
$$
Therefore, setting $\delta=\frac{\varphi_2-\varphi_1}{2}$ we have
\begin{equation}\label{complex}
r_1e^{-i\delta}= r(\bar{\varphi})+a(\dot{r}(\bar{\varphi})
+ir(\bar{\varphi})),\quad
r_2e^{i\delta}=r(\bar{\varphi})+b(\dot{r}(\bar{\varphi})
+ir(\bar{\varphi})).
\end{equation}
From these two equations we find
$$a=-\frac{ r(\bar{\varphi})\tan\delta}{r(\bar{\varphi})+\dot{r}(\bar{\varphi})\tan\delta},\qquad
b=\frac{ r(\bar{\varphi})\tan\delta}{r(\bar{\varphi})-\dot{r}(\bar{\varphi})\tan\delta}.
$$
and then equating imaginary parts of (\ref{complex}) we obtain the
formulas:
\begin{equation}
\label{r12}
r_1=\frac{r^2(\bar{\varphi})}{r(\bar{\varphi})\cos\delta+\dot{r}(\bar{\varphi})\sin
\delta},\qquad
r_2=\frac{r^2(\bar{\varphi})}{r(\bar{\varphi})\cos\delta-\dot{r}(\bar{\varphi})\sin
\delta}.
\end{equation}
 Formulas (\ref{r12}) lead us to the following

\vskip5mm

\noindent {\bf Theorem 5} {\it The Angular billiard map $\ca$
preserves the symplectic form $d\left(\frac{1}{r}\right)\wedge
d\varphi$ and can be written with the help of the generating
function
$$S(\varphi_1,\varphi_2)=\frac{2}{r(\bar{\varphi})}
\sin\delta=\frac{2}{r(\frac{\varphi_1+\varphi_2}{2})}
\sin(\frac{\varphi_2-\varphi_1}{2}),$$  $$ -\frac{\partial
S}{\partial {\varphi_1}}=\frac{1}{r_1},\qquad \frac{\partial
S}{\partial {\varphi_2}}=\frac{1}{r_2}.
$$

Moreover the Angular billiard map is a twist map near the boundary
$\Gamma,$ i.e. function $S$ satisfies the twist condition:
$$\frac{\partial^2S}{\partial\varphi_1\partial\varphi_2}>0.$$
}

\begin{proof}
Let $p=\frac{1}{r},\ p(\varphi)=\frac{1}{r(\varphi)}.$
Differentiating the function $S$ and comparing with (\ref{r12}) we
have
\begin{equation}
-\frac{\partial S}{\partial {\varphi_1}}=p(\bar{\varphi})\cos
\delta-\dot{p}(\bar{\varphi})\sin\delta=\frac{1}{r_1},\qquad
\frac{\partial S}{\partial {\varphi_2}}=p(\bar{\varphi})\cos
\delta+\dot{p}(\bar{\varphi})\sin\delta=\frac{1}{r_2}.
\end{equation}
But this exactly means that $S$ is the generating function. Next we
compute:
$$
\frac{\partial^2S}{\partial\varphi_1\partial\varphi_2}(\varphi_1,\varphi_2)=\frac{1}{2}(\ddot{p}(\bar\varphi)+p(\bar\varphi))\sin\delta.
$$
But the expression in brackets is positive since it equals the
curvature radius of $\gamma$ in terms of supporting function.

\end{proof}

\subsection{Reduction to KAM for angular billiards}

 Let us
remark first that the invariant curves of  Angular billiard are dual
to the caustics of the Birkhoff billiard. The proof of existence of
invariant curve for angular billiard or equivalently existence of
caustics for inner billiard is traditionally reduced, after
Lazutkin, to the following KAM result  (see \cite{L}). Usually
people use the second choice of symplectic coordinates for this
reduction and this leads to very complicated formulas (see also
\cite{siburg}). Our new observation here is that in the coordinates
$(\varphi,p)$ the reduction to the KAM lemma of Lazutkin
 looks much simpler.

\vskip5mm

\noindent{\bf Theorem} (Lazutkin, 1972) {\it Given a transformation
$T$ of the cylinder
written in coordinates $(x\  {\rm mod}\ 1,\ y\geq0)$ in the form
$$
 x'=x+y+y^{m+1}f(x,y),\quad y'=y+y^{m+2}g(x,y),\quad m\geq 1,
$$
with sufficiently smooth function $f,g.$ Assume that for any curve
$\alpha$ homotopic to the curve $\{y=0\}$ and sufficiently close to
it, $T(\alpha)\cap\alpha\neq\emptyset.$ Then $T$ has invariant
curves arbitrary close to ${y=0}.$ }

\vskip5mm

 Application of this  theorem can be
done via the following Lemma \cite{T1} (the intersection property is
obviously satisfied for billiards, by invariance of symplectic
form):

\vskip5mm

\noindent{\bf Lemma} (see \cite{T1}) {\it Let $T$ be a
transformation of the cylinder which in coordinates $(u\ {\rm mod} \
1,\ v\geq 0)$ has the form
$$
 u'=u+\alpha(u)\ v+(v^2),\quad v'=v+\beta(u)\ v^2+(v^3),
$$
where the function $\alpha$ has no zeros. Then their exists change of
coordinates
$$
x(u,v)=a(u)+(v),\quad y(u,v)=h(u)+(v^2),
$$
in which $T$ has the form as in the Theorem:
$$x'=x+y+(y^2),\quad y'=y+(y^3).$$
Here $(f)$ denotes a smooth function in the ideal generated by $f.$
}
\vskip5mm

 In order to apply this Lemma to
angular billiard map $\ca$ we work in coordinates $(p,\varphi),$
where $p=\frac{1}{r}.$ First we change coordinates from $(p,\varphi)$ to $(z,\varphi),$ where $z=p(\varphi)-p,$ so that $z\geq 0$
exactly in the complement of the curve $\Gamma.$ By formulas
(\ref{r12}) transformation $\ca:(z_1,\varphi_1)\rightarrow
(z_2,\varphi_2)$ acts according to implicit formulas
\begin{equation}\label{z1}
z_1=p(\varphi_1)-p(\bar{\varphi})\cos
\delta+\dot{p}(\bar{\varphi})\sin\delta,\quad z_2=p(\varphi_2)-
p(\bar{\varphi})\cos \delta-\dot{p}(\bar{\varphi})\sin\delta.
\end{equation}
Rename now in (\ref{z1})
$$
 \varphi_1\rightarrow \varphi;\quad \bar{\varphi}\rightarrow \varphi+\delta;
 \quad \varphi_2\rightarrow\varphi'=
 \varphi+2\delta;\quad z_1\rightarrow z,\quad z_2\rightarrow z'.
$$
Then (\ref{z1}) can be written as follows:
\begin{equation}\label{z2}
z=p(\varphi)-p(\varphi+\delta)\cos\delta+\dot{p}(\varphi+\delta)\sin\delta,
\end{equation}
\begin{equation}\label{z3}
z'=z+p(\varphi+2\delta)-p(\varphi)-2\dot{p}(\varphi+\delta)\sin\delta.
\end{equation}
From (\ref{z2}) we get an expansion in power series in $\delta$:
\begin{equation}
\label{expand}
z=A\delta^2+B\delta^3+(\delta^4),\quad
for\quad  A=\frac{1}{2}(\ddot{p}(\varphi)+p(\varphi)),
\quad B=\frac{2}{3}\dot{A}.
\end{equation}
In the sequel it is important that $A$ is positive everywhere.

Put $v=\sqrt z$ and write $\delta=a\cdot v+(v^2)$ and in order to find $a$ we substitute back into (\ref{expand}) and find out that
$$
 a=\frac{1}{\sqrt A}.
$$
Thus we have
$$
\delta=\frac{1}{\sqrt A}v+(v^2).
$$
Let us recall that $\delta=\frac{\varphi'-\varphi}{2}$ and thus
\begin{equation}\label{star1}
\varphi'=\varphi+\frac{2}{\sqrt A}v+(v^2).
\end{equation}
Take now equation (\ref{z3}), we have
$$(v')^2=v^2+B\cdot\delta^3+(\delta^4)=v^2+B(\frac{1}{\sqrt A}v+(v^2))^3+(v^4)=v^2+\frac{B}{A^{3/2}}v^3+(v^4).$$
Extracting square root of the last equation, we get
\begin{equation}
\label{star2}
v'=v+\frac{B}{2A^{3/2}}v^2+(v^3).
\end{equation}
So we have the formulas (\ref{star1}) and (\ref{star2})  which are
precisely as in the lemma. Therefore KAM result by Lazutkin applies
giving the existence of invariant curves and we are done.

\section{Polynomial integrals}
Here we recall known results from \cite{B}, \cite{KT}.

Let $\Phi(x,v), x\in\Omega,\ |v|=1,$ be a non-constant function on the
unite tangent bundle $T_1\Omega,$ where $\Omega$ is the Birkhoff
billiard table bounded by the curve $\gamma$ in the plane. Assume
that $\Phi$ is an integral of the billiard flow on the energy level
$|v|=1$ which is a polynomial of degree $n$ in the components of
vector $v.$ Then it follows from \cite{B} and \cite{KT}
that there exists a non-constant \emph{homogeneous} polynomial which is
also an integral of the billiard flow. This is done as follows:
split first $\Phi$ into two parts of all even
degrees and all odd degrees, $\Phi=\Phi_{even}+\Phi_{odd}.$ Using the involution
$\iota:(x,v)\mapsto(x,-v),$ we have that the function
$\Phi\circ\iota=\Phi_{even}-\Phi_{odd}$ is also an integral of the
billiard flow. Therefore, $\Phi_{even}$ and $\Phi_{odd}$ separately
are also integrals. The last thing to do is to homogenize
$\Phi_{even}$ and $\Phi_{odd}$ by multiplying their summands by
$v^2=v_x^2+v_y^2$ in appropriate power.

Therefore, we may assume
that $\Phi(q,v)$ is a homogeneous polynomial in $v_x,v_y.$ It then
follows that $\Phi$ is in fact homogeneous polynomial in
$v_x,v_y,\sigma$ where $\sigma=v_yq_1-v_xq_2$ is a momentum of the
vector $v$ about the origin. Indeed, since
$$
 \Phi(q+tv,v)=\Phi(q,v)
$$
we have
$$
 \frac{d}{dt}\Phi(q+tv,v)=v_x\Phi_{q_1}(q+tv,v)+v_y\Phi_{q_2}(q+tv,v)=0,
$$
so
$$
 v_x\Phi_{q_1}(q,v)+v_y\Phi_{q_2}(q,v)=0.
$$
 Hence,
$$
 \Phi(q,v)=\Phi(\sigma,v_x,v_y),\qquad \sigma=q_2v_x-q_1v_y
$$
and $\Phi(\sigma,v_x,v_y)$ is a homogeneous polynomial in
$\sigma,v_x,v_y$ of degree $n.$

Let $\tau(q)$ is a unit tangent vector to $\gamma$ at $q\in\gamma$
with positive orientation of $\gamma$ then it is easy to see that
\begin{equation}
\label{tau}
 \Phi(q,\tau(q))=c,\ \forall q\in\gamma,\qquad c\in{\mathbb R},
\end{equation}
where $c$ is a constant. One can see this considering a sequence of
trajectories which tend to the boundary.

Moreover, for even $n$ we can assume that the constant $c$ equals 0.
Indeed, instead of $\Phi(q,v)$ we can consider the integral
$\Phi(q,v)-c|v|^n$ which is also a polynomial integral for even $n.$
And if $n$ is odd then one can consider the integral $\Phi^2$  and
proceed as above. So, we shall assume from now on that $n$ is even
and
\begin{equation}\label{i1}
 \Phi(q,\tau(q))=0,\qquad q\in\gamma.
\end{equation}

 It then follows from Theorem 3 that the dual  $\Gamma$ is real oval of an
algebraic curve, consequently, the dual curve to $\Gamma,$ i.e.
$\gamma,$ is also a real oval of an algebraic curve.

\vskip5mm

%%%%%%%%%%%%%%%%%%%%%%%%%%%%%%%%
\noindent $\mathbf{Remark.}$ {\it Let Birkhoff billiard flow inside
$\gamma$ admits a polynomial integral $\Phi$ of the form
$$\Phi(\sigma,v)=\sum_{k+l+m=n}f_{klm}(v_x)^k(v_y)^l\sigma^m.$$
Then the function
$$
\hat{\Phi}(\varphi,
p)=\sum_{k+l+m=n}f_{klm}(-\sin\varphi)^k(\cos\varphi)^lp^m,
$$ is invariant under Birkhoff billiard mapping $\cb:(\varphi,p)\mapsto(\varphi',p'),$
since
$$
 \sigma=p|v|,\quad v_x=-|v|\sin\varphi,\quad v_y=|v|\cos\varphi.
$$
Notice that $\hat{\Phi}$ is polynomial function in $p$ with
coefficients which are trigonometric polynomials of $\varphi.$ In
fact this correspondence is valid also in the opposite direction.
Given an integral of the mapping $\cb$ which is polynomial in $p$,
   whose coefficients are trigonometric polynomials, one can reconstruct polynomial
   integral of the billiard flow.
   Thus these two problems are equivalent. We refer to  survey
    \cite{V} on integrable mappings.

Let us stress that it is a much harder problem to study first
integrals of $\cb$ which are polynomial in $p$ with smooth periodic
coefficients (which are not necessarily trigonometric polynomials).}

\section{Proof of Theorem 3}

Given a homogenous polynomial integral $\Phi$ of even degree for
Birkhoff billiard flow. We define the function $G$ at the vicinity
of the dual curve $\Gamma$ as follows. Take a point $A\in
U\setminus\cs,$ let $a$ be the dual line to $A$ and let $v_a$ be a
unite tangent vector to $a$ with a positive momentum
$\sigma(v_a)=p>0.$ Then set
$$G(A):=\Phi(\sigma(v_a),v_a).$$
Notice that since $\Phi$ is a first integral then $G$ is correctly
defined, that is it does not depend on the foot point for $v_a$ on
the line $a.$ Moreover, let $B=\ca (A)$ then it follows from Theorem
2 that for the dual line $b$ there are two possibilities:
$$v_b=\pm\cb(v_a).$$
But since $\Phi$ is assumed to be of even degree we have for
every $v:$
$$
\Phi(\sigma(v),v))=\Phi(\sigma(-v),-v).
$$
Therefore,
$$G(B)=\Phi(\sigma(v_b),v_b)=\Phi(\pm\cb (v_a))=\Phi(\sigma(v_a),v_a)=G(A).$$
So $G$ is an integral of Angular billiard map $\cb.$ Moreover, since
by (\ref{i1}), $\Phi$ vanishes on tangent vectors to $\gamma,$ then
$G$ vanishes on $\Gamma$ by the definition.

Now let us compute $G$ in coordinates.
 Let $v$ be a unite tangent vector along the line $a$ with positive momentum
 $\sigma(v)=p>0.$ We have:
$$\Phi(\sigma(v),v)=\sum_{k+l+m=n}f_{klm}(v_x)^k(v_y)^lp^m.$$
Recall that the dual point $A$ of the line $a$ has the coordinates,
see (\ref{duality}):
$$x=\frac{n_1}{p}=\frac{\cos\varphi}{p}=\frac{v_y}{p},\quad y=\frac{n_2}{p}
=\frac{\sin\varphi}{p}=\frac{-v_x}{p}.$$
In particular
$$p=\frac{1}{r}=\frac{1}{\sqrt{x^2+y^2}}.$$
Then we compute:
$$
 G(x,y)=\sum_{k+l+m=n}f_{klm}(-py)^k(px)^lp^m=
 p^n\left(\sum_{k+l+m=n}f_{klm}(-y)^k(x)^l\right)=
$$
$$
 =p^nF(x,y)=
 \frac{F(x,y)}{(x^2+y^2)^{n/2}},
$$
where the polynomial in brackets is
denoted by $F.$ Since $G$ vanishes on $\Gamma$ then also $F$ does.
This is exactly the form which was required in Theorem 3. This
completes the proof.

\vskip5mm \noindent$\mathbf {Remark}.$ {\it If the degree $n$ of the
polynomial integral $\Phi$ is odd one still can define $G$ as above.
But in this case $G$ is not an integral of $\ca$ globally,
but only near the boundary $\Gamma$.

However, one cannot claim anymore that $G$
necessarily vanishes on the boundary $\Gamma,$ but only that it
equals to a constant there, by (\ref{tau}).}

\section{Proof of Theorem 1}
In this Section we derive a remarkable identity and then study it in
algebro-geometric terms.

 Let $\Gamma$ be a dual curve of the curve
$\gamma,$ and assume that Birkhoff billiard inside $\gamma$ admits a
polynomial integral of an even degree $2p$, where $p$ is an integer.
Then according to Theorem 2 Angular billiard for $\Gamma$ is
integrable with the integral of the form
$$
 G_1(x,y)=\frac{F_1(x,y)}{(x^2+y^2)^p}.
$$
Here $F_1$ is a polynomial of $\deg F_1=2p,$ and $F_1$ is assumed to
vanish on $\Gamma.$

Let $\Gamma$ is defined by the equation $f=0,$ where $f$ is an
irreducible polynomial in $\mathbb{C}[x,y]$ of degree $d.$
  Since $F_1=0$ on
$\Gamma$ one can write $F_1$ in the form:
$$
 F_1(x,y)=f^k(x,y)g_1(x,y),
$$
where $k$ is positive integer, and $g_1\ne 0$ identically on
$\Gamma.$  It is important, that $f,g_1$ can be assumed to be
\emph{real} polynomials. Let $\Gamma_1$ be the arc of $\Gamma$ where
$g_1> 0$ (if $g_1<0$ everywhere we change the sign of $G_1$). Let
$$
 F(x,y)=(F_1(x,y))^{\frac{1}{k}}=f(x,y)g(x,y),\qquad g(x,y)=(g_1(x,y))^{\frac{1}{k}}.
$$
Next we replace $G_1$ by $G:=G_1^{\frac{1}{k}}$:
$$
 G(x,y)=\frac{(F_1(x,y))^{\frac{1}{k}}}{(x^2+y^2)^{\frac{p}{k}}}:
 =\frac{F(x,y)}{(x^2+y^2)^m},\qquad m=\frac{p}{k}.
$$
Then $G$ is also an integral, which also vanishes on $\Gamma_1,$ but
$F,g$ are not necessarily polynomials anymore.

\vskip8mm

\begin{lemma}
\label{automs}
The identity
\begin{equation}\label{e1}
 F(x+\varepsilon F_y,y-\varepsilon F_x) \left(-\frac{\mu}{\varepsilon}\right)^{2m}=
 F(x+\mu F_y,y-\mu F_x)
\end{equation}
holds for small real $\varepsilon,$ where $(x,y)\in\Gamma_1$ and
$$
 \mu=-\frac{(x^2+y^2)\varepsilon}{x^2+y^2+2\varepsilon(xF_y-yF_x)}.
$$
\end{lemma}

\vskip8mm

\begin{proof}
The vector
$$v=(F_y,-F_x)=(g_yf+gf_y,-g_xf-gf_x)=(gf_y,-gf_x)$$
is tangent vector to  $\Gamma_1$ at the point $T(x,y).$ Let $A,B$ be
two points such that $\overrightarrow{OA}=\varepsilon
v,\overrightarrow{OB}=\mu v.$ Assume that $\ca(A)=B$ so the
angle $\alpha=\angle AOT$ is equal to $\beta=\angle TOB,$ where
$O=(0,0).$ Then we have
$$
 \frac{(\overrightarrow{OA},\overrightarrow{OT})^2}{\vert OA\vert^2\vert OT\vert^2}=\frac{(\overrightarrow{OB},\overrightarrow{OT})^2}{\vert OB\vert^2\vert OT\vert^2},
$$
where
$$
 \overrightarrow{OA}=(x+\varepsilon F_y,y-\varepsilon F_x),\qquad \overrightarrow{OB}=(x+\mu F_y,y-\mu F_x),
$$
$$
 (\overrightarrow{OA},\overrightarrow{OT})=x^2+y^2+\varepsilon (xF_y-yF_x),\qquad (\overrightarrow{OB},\overrightarrow{OT})=x^2+y^2+\mu (xF_y-yF_x),
$$
$$
 \vert OA\vert^2=x^2+y^2+2\varepsilon(xF_y-yF_x)+\varepsilon^2(F_x^2+F_y^2),\qquad
 \vert OB\vert^2=x^2+y^2+2\mu(xF_y-yF_x)+\mu^2(F_x^2+F_y^2).
$$
So, we obtain the equation on $\mu$
$$
 \frac{(x^2+y^2+\varepsilon (xF_y-yF_x))^2}{(x^2+y^2+2\varepsilon(xF_y-yF_x)+
 \varepsilon^2(F_x^2+F_y^2))^2}= \frac{(x^2+y^2+\mu (xF_y-yF_x))^2}{(x^2+y^2+2\mu(xF_y-yF_x)+\mu^2(F_x^2+F_y^2))^2}.
$$
This equation has two solutions: the trivial one: $\mu=\varepsilon,$
and the relevant one for us:
$$
 \mu=-\frac{(x^2+y^2)\varepsilon}{x^2+y^2+2\varepsilon(xF_y-yF_x)}.
$$
The condition that $G$ is an integral for $\ca$ reads:
$$
 G(x+\varepsilon F_y,y-\varepsilon F_x)=G(x+\mu F_y,y-\mu F_x),
$$
or equivalently

$$
 \frac{F(x+\varepsilon F_y,y-\varepsilon F_x)}{((x+\varepsilon F_y)^2+(y-\varepsilon F_x)^2)^m}
=\frac{F(x+\mu F_y,y-\mu F_x)}{((x+\mu F_y)^2+(y-\mu F_x)^2)^m}.
$$
Consider now the triangle $\triangle AOB$ where $OT$ is the bisector. Then by known property of bisector
we have
$$\frac{|OA|}{|OB|}=\frac{|TA|}{|TB|}.
$$
Therefore
$$
 \sqrt{\frac{(x+\varepsilon F_y)^2+(y-\varepsilon F_x)^2}{(x+\mu F_y)^2+(y-\mu F_x)^2}
 }=-\frac{\varepsilon}{\mu}.
$$
This implies the formula (\ref{e1}) of the lemma.
\end{proof}

We shall use the following notation. For any function $f$ we set
$$H(f):=f_y(f_{xx}f_y-f_{xy}f_x)+f_x(f_{yy}f_x-f_{xy}f_y)=f_{xx}f_y^2-2f_{xy}f_xf_y+f_{yy}f_x^2
.$$
It turns out that the identity (\ref{e1}) can be integrated to give
the following remarkable formula:

\vskip8mm

\begin{theo}
\label{remarkable} The following formula holds true for all $(x,y)\in\Gamma_1:$
\begin{equation}\label{e2}
 g^3(x,y)H(f(x,y))=c_1(x^2+y^2)^{3m-3},
\end{equation}
 where $c_1$ is a constant.
\end{theo}

\vskip8mm

\begin{proof}
First notice that the expression for $\mu$ in Lemma \ref{automs}
yields the power series for $\mu$:
$$\mu=\sum_{k=1}^{\infty}\mu_k\veps^k,\quad \mu_1=-1,\quad
\mu_k=(-1)^k\left(2\frac{xF_y-yF_x}{x^2+y^2}\right)^{k-1}.
$$
Substituting the series into (\ref{e1}) we get the following:
\begin{equation}
\label{e7} F(x+\varepsilon F_y,y-\varepsilon
F_x)(1-\mu_2\veps\dots)^{2m}- F(x+(-\veps+\mu_2\veps^2\dots)
F_y,y-(-\veps+\mu_2\veps^2\dots) F_x)=0.
\end{equation}

It is easy to see that the coefficients at
$\varepsilon^0,\varepsilon,\varepsilon^2$ in the left hand side of
(\ref{e7}) vanish for all $(x,y)\in\Gamma_1.$ Equating to zero the
coefficient at $\veps^3$ of (\ref{e7}) we get the equation:
$$\frac{1}{3}(F_{xxx}F_y^3-3F_{xxy}F_y^2F_x+3F_{xyy}F_yF_x^2-F_{yyy}F_x^3)+
\frac{1}{2}(F_{xx}F_y^2-2F_{xy}F_xF_y+F_{yy}F_x^2)(-2m\mu_2)-
$$
$$
-\frac{1}{2}(F_{xx}F_y^2-2F_{xy}F_xF_y+F_{yy}F_x^2)(-2\mu_2)=0.
$$
Substituting $\mu_2$ into the last equation we have:
$$
(F_{xxx}F_y^3-3F_{xxy}F_y^2F_x+3F_{xyy}F_yF_x^2-F_{yyy}F_x^3)+
$$
\begin{equation}
\label{e8}
\\
+6(F_{xx}F_y^2-2F_{xy}F_xF_y+F_{yy}F_x^2)(1-m)\left(\frac{xF_y-yF_x}{x^2+y^2}\right)=0.
\end{equation}
Denote by $v=F_y\partial_x-F_x\partial_y$ be the vector field
tangent to $\Gamma_1$ at $(x,y).$ Recall that
$$
H(F)=F_{xx}F_y^2-2F_{xy}F_xF_y+F_{yy}F_x^2,
$$
and notice that the following two identities hold along $\Gamma_1:$
$$L_v(x^2+y^2)=2(xF_y-yF_x),$$
$$
L_vH(F)=F_{xxx}F_y^3-3F_{xxy}F_y^2F_x+3F_{xyy}F_yF_x^2-F_{yyy}F_x^3.
$$
Using these identities equation (\ref{e8}) can be rewritten as follows:
\begin{equation}\label{e9}
(x^2+y^2)L_vH(F)+3H(F)(1-m)L_v(x^2+y^2)=0.
\end{equation}
Multiplying (\ref{e9}) by $(x^2+y^2)^{2-3m}$ we get
\begin{equation}
\label{e10}
(x^2+y^2)^{3-3m}L_vH(F)+3H(F)(1-m)(x^2+y^2)^{2-3m}L_v(x^2+y^2)=0.
\end{equation}
But the left hand side of (\ref{e10}) is the complete derivative and
thus
$$L_v\left(H(F)(x^2+y^2)^{3-3m}\right)=0.$$

Since $v$ is a tangent vector field to
$\Gamma_1,$ then the function $H(F)(x^2+y^2)^{3-3m}$ must
be a constant on $\Gamma_1,$ i.e.
\begin{equation}\label{eq2}
 H(F(x,y))=c_1(x^2+y^2)^{3m-3},\qquad \forall(x,y)\in\Gamma_1,
\end{equation}
for some constant $c_1.$ Since $F=f g$ and, in addition, one can check
by a direct calculation that
$$
 H(f(x,y)g(x,y))=g^3(x,y)H(f(x,y)),
$$
we conclude with the formula
$$
g^3(x,y)H(f(x,y))=c_1(x^2+y^2)^{3m-3},\qquad \forall(x,y)\in\Gamma_1.
$$
This proves Theorem \ref{remarkable}.
\end{proof}

 We shall use the following notations. For any polynomial  $p(x,y),$
we denote by $\tilde{p}(x,y,z)$ the corresponding homogeneous
polynomial of the same degree as $p.$

\vskip8mm

\begin{lemma}
The identity
\begin{equation}\label{e4}
 \tilde{g}_1^6(x,y,z)({\rm Hess}(\tilde{f}(x,y,z)))^{2k}+
 c(x^2+y^2)^{6p-6k}=\tilde{f}(x,y,z)\tilde{h}(x,y,z),
\end{equation}
holds true for all $(x,y,z)\in\mathbb{C}^3,$ where $c$ is a
constant, $\tilde{h}$ is a homogeneous polynomial, and
$$
 {\rm Hess}(\tilde{f}(x,y,z)):=\det\left(
  \begin{array}{ccc}
    \tilde{f}_{xx} & \tilde{f}_{xy} & \tilde{f}_{xz} \\
    \tilde{f}_{xy} & \tilde{f}_{yy} & \tilde{f}_{yz} \\
    \tilde{f}_{xz} & \tilde{f}_{yz} & \tilde{f}_{zz} \\
  \end{array}
\right).
$$
\end{lemma}

\vskip8mm
\begin{proof}
Let us raise both sides of (\ref{e2}) to the power $2k.$ This gives
$$
 g^{6k}(x,y)(H(f(x,y)))^{2k}=g_1^{6}(x,y)(H(f(x,y)))^{2k}=
 c_1^{2k}(x^2+y^2)^{6mk-6k}
$$
or
$$
 g_1^{6}(x,y)(H(f(x,y)))^{2k}=
 c_1^{2k}(x^2+y^2)^{6p-6k},\qquad \forall(x,y)
 \in\Gamma_1.
$$
Since on both sides of the previous identity we have polynomials
then their difference is divisible by $f.$ Therefore, we get
\begin{equation}\label{e5}
  g_1^6(x,y)(H(f(x,y)))^{2k}-c_1^{2k}(x^2+y^2)^{6p-6k}=f(x,y)h_1(x,y),\qquad (x,y)\in{\mathbb C}^2,
\end{equation}
where $h_1(x,y)$ is a polynomial. In order to homogenize, let us
note that
$$
 \deg  g_1^{6}(x,y)(H(f(x,y)))^{2k}=\deg (x^2+y^2)^{6p-6k}+4k.
$$
Indeed,
$$\deg (x^2+y^2)^{6p-6k}=12p-12k.
$$
Furthermore, denoting $\deg f=d$ and $\deg g_1=q,$ we have
$$
 \deg F_1=\deg f^k+\deg g_1=2p=kd+q.
$$
Hence,
$$
 \deg  g_1^{6}(H(f))^{2k}=6q+2k(3d-4)=6q+6kd-8k=12p-8k.
$$
From (\ref{e5}) we obtain
\begin{equation}\label{e6}
 \tilde{g}_1^6(x,y,z)(H(\tilde{f}(x,y,z)))^{2k}-c_1^{2k}z^{4k}(x^2+y^2)^{6p-6k}=
 \tilde{f}(x,y,z)\tilde{h}_1(x,y,z).
\end{equation}
Now we use the identities (see \cite{W}):
$$
 {\rm Hess}(\tilde{f})=\frac{(d-1)^2}{z^2}
  \det\left(
  \begin{array}{ccc}
    \tilde{f}_{xx} & \tilde{f}_{xy} & \tilde{f}_{x} \\
    \tilde{f}_{xy} & \tilde{f}_{yy} & \tilde{f}_{y} \\
    \tilde{f}_{x} & \tilde{f}_{y} & \frac{d}{d-1}\tilde{f}\\
  \end{array}
\right)=$$
\begin{equation}\label{Hf}
=\frac{(d-1)^2}{z^2}\left(\frac{d}{d-1}\tilde{f}
 (\tilde{f}_{xx}\tilde{f}_{yy}-\tilde{f}_{xy}^2)-H(\tilde{f})\right)
\end{equation}
where $H(\tilde{f}(x,y,z))$ at the right hand side is computed with
$z$ being a parameter. The last step is just to substitute the
expression for $H(\tilde{f})$ via ${\rm Hess}(\tilde{f})$ from
(\ref{Hf}) into (\ref{e6}) in order to get the required (\ref{e4}).
\end{proof}

 Now we are ready to finish the
proof of Theorem 1. We follow the idea of Lemma 3 of \cite{T}.
Consider the situation in $\mathbb{C}P^2.$ Any intersection point in $\mathbb{C}P^2$ between Hessian curve of
${\rm Hess}(\tilde{\Gamma})$ with
$\tilde{\Gamma}$ is either singular or inflection point of $\tilde{\Gamma}$. So, if there is
a singular or inflection point $(x_0:y_0:z_0)\in\tilde{\Gamma}$ such
that $x_0^2+y_0^2\ne 0,$ it then follows from (\ref{e4}) that $c=0.$
Therefore, ${\rm Hess}(\tilde{f})$  must vanish on $\tilde{\Gamma},$ since
$\tilde{g}_1\ne 0$ identically on $\tilde{\Gamma}.$ But this can
happen is only if $\tilde{\Gamma}$ is a line (see \cite{F}), but this is impossible.

Let us prove now that $\tilde{\Gamma}$ must have
singular points. If on the contrary $\tilde{\Gamma}$ is a smooth
curve, then it follows from (\ref{e4}) that all inflection points
must belong to two lines $L_1$ and $L_2$ defined by the equations
$$
 L_1=\{ x+iy=0\},\qquad\ L_2=\{x-iy=0\}.
$$
Recall, $d$ is the degree of $\tilde{\Gamma}.$ Then the Hessian
curve intersects $\tilde{\Gamma}$ exactly in inflection points, and
moreover, it is remarkable fact that the intersection multiplicity
of such a point of intersection equals exactly the order of
inflection point (see \cite{W}), and hence does not exceed  $d-2.$
Furthermore, the lines $L_1$ and $L_2$ intersect $\tilde{\Gamma}$
maximum in $2d$ points together . Hence,  we have altogether counted
with multiplicities not more than $ 2d(d-2),$ but on the other hand
the Hessian curve has degree $3(d-2)$ and thus by Bezout theorem the
number of intersection points with multiplicities is $3d(d-2).$ This
contradiction shows that $\tilde{\Gamma}$ can not be a smooth curve
unless $d=2.$ Theorem 3 is proved.

\vskip8mm

\noindent M. Bialy, {\it email:} bialy@post.tau.ac.il
\vskip5mm
\noindent A.E. Mironov, {\it email:} mironov@math.nsc.ru

\end{document}